\documentclass[letterpaper,10pt,conference,draftcls,onecolumn]{ieeeconf}
\IEEEoverridecommandlockouts
\usepackage{cite}
\usepackage{amsmath,amssymb,amsfonts}
\usepackage{algorithmic,algorithm}
\usepackage{graphicx}
\usepackage{textcomp}
\usepackage{xcolor}
\usepackage{enumerate}
\usepackage{mathtools}

\newcommand{\R}{\mathbb{R}}
\newcommand{\N}{\mathbb{N}}

\newtheorem{assumption}{Assumption}
\newtheorem{theorem}{Theorem}
\newtheorem{problem}{Problem}
\newtheorem{lemma}{Lemma}
\newtheorem{definition}{Definition}
\newtheorem{remark}{Remark}

\DeclareMathOperator*{\argmax}{arg\,max}
\DeclareMathOperator*{\argmin}{arg\,min}

\newcommand{\mcc}{\mu^c_c}
\newcommand{\mhat}{\hat{\mu}_{c, \delta}}
\newcommand{\xhat}{\hat{x}_{\delta}}
\newcommand{\xhatt}{\hat{x}^t}

\newcommand{\mhatd}{\hat{\mu}_{\delta}}

\usepackage{mathtools}
\mathtoolsset{showonlyrefs=true}

\usepackage{cite}

\usepackage{color}

\usepackage{marginnote}

\def\BibTeX{{\rm B\kern-.05em{\sc i\kern-.025em b}\kern-.08em
    T\kern-.1667em\lower.7ex\hbox{E}\kern-.125emX}}
\begin{document}

\title{Towards Totally Asynchronous Primal-Dual Convex Optimization in Blocks\\
\thanks{This work was supported by a task order contract from the Air Force Research Laboratory through
Eglin AFB.}}

\author{
Katherine R. Hendrickson and Matthew T. Hale$^*$\thanks{${}^*$The authors are with the Department of Mechanical and Aerospace Engineering,
Herbert Wertheim College of Engineering, University of Florida. Emails:
\texttt{\{kat.hendrickson,matthewhale\}@ufl.edu}. 
}
}

\maketitle

\begin{abstract}
We present a parallelized primal-dual algorithm for solving constrained convex optimization problems. The algorithm is ``block-based,'' in that vectors of primal and dual variables are partitioned into blocks, each of which is updated only by a single processor. We consider four possible forms of asynchrony: in updates to primal variables, updates to dual variables, communications of primal variables, and communications of dual variables. We explicitly construct a family of counterexamples to rule out permitting asynchronous communication of dual variables, though the other forms of asynchrony are permitted, all without requiring bounds on delays. A first-order update law is developed and shown to be robust to asynchrony. We then derive convergence rates to a Lagrangian saddle point in terms of the operations agents execute, without specifying any timing or pattern with which they must 
be executed. These convergence rates contain a synchronous algorithm as a special case and are used to quantify an ``asynchrony penalty.'' Numerical results illustrate these developments. 
\end{abstract}

\section{Introduction}
A wide variety of machine learning problems can be formalized
as convex programs~\cite{sra12,shalev12,bertsekas15,boyd04}. Large-scale machine learning then requires
solutions to large-scale convex programs, which
can be accelerated through parallelized solvers running
on networks of processors. 
In large networks, it is difficult (or
outright impossible) to synchronize 
their behaviors.
The behaviors of interest are computations, which
generate new information, and communications, which
share this new information with other processors. Accordingly,
we are interested in asynchrony-tolerant large-scale optimization.

The challenge of asynchrony is that
it causes disagreements
among processors that result from receiving different
information at different times. One way to reduce
disagreements
is through repeated averaging of processors' iterates.
This approach dates back several decades~\cite{tsitsiklis86}, and
approaches of this class include primal~\cite{nedic09,nedic10,nedic15},
dual~\cite{duchi11,tsianos12,tsianos12b}, and primal-dual~\cite{zhu11,jaggi14} algorithms. 
However, these averaging-based methods require bounded
delays in some form, often through requiring connectedness
of agents' 
communication 
graphs over intervals of a prescribed length~\cite[Chapter 7]{bertsekas89}.
In some applications, delays are outside agents' control, e.g.,
in a contested environment where communications are jammed,
and thus delay bounds cannot be reliably enforced.
Moreover, graph connectivity cannot be easily checked
locally by individual agents, meaning even satisfaction
or violation of connectivity bounds is not readily ascertained. 
In addition, these methods require multiple processors to 
update each decision variable, which duplicates computations
and increases processors' workloads.
This can be prohibitive in large problems, such as learning 
problems with billions of data points. 

Therefore, in this paper we develop a parallelized primal-dual
method for solving large constrained convex optimization problems. 
Here, by ``parallelized,'' we mean that each decision variable
is updated only by a single processor. As problems grow, this
has the advantage of keeping each processor's computational
burden approximately constant. The decision variables assigned
to each processor are referred to as a ``block,'' and
block-based algorithms date back several decades as well~\cite{tsitsiklis86,Bertsekas1983}.
Those early works solve unconstrained or set-constrained
problems, in addition to select
problems with functional constraints~\cite{bertsekas89}. 
To bring parallelization to arbitrary constrained
problems, we develop a primal-dual approach that does not
require constraints to have a specific form.

Block-based methods have previously been shown to tolerate
arbitrarily long delays in both communications and
computations in unconstrained problems~\cite{Bertsekas1983,hochhaus18,ubl19}, 
eliminating the need to enforce
and verify delay boundedness assumptions. For constrained
problems of a general form, block-based methods have been
paired with primal-dual algorithms with centralized
dual updates~\cite{hale14,hale17} and/or synchronous primal-updates~\cite{koshal2011multiuser}.
To the best of our knowledge, 
arbitrarily asynchronous block-based updates have not been
developed for convex programs of a general form. 
A counterexample in~\cite{hale17} suggested that arbitrarily asynchronous
communications of dual variables can preclude convergence, though
that example leaves open the extent to which dual asynchrony
is compatible with convergence.

In this paper, we present 
a primal-dual optimization algorithm that permits
arbitrary asynchrony in primal variables, while
accommodating dual asynchrony to the extent possible.
Four types of asynchrony are possible: (i) asynchrony
in primal computations, 
(ii) asynchrony in communicating primal variables,
(iii) asynchrony in dual computations,
(iv) asynchrony in communicating dual variables. 
The first contribution of this paper is to show that 
item (iv) is fundamentally problematic 
 using an explicit family of counterexamples that we construct. 
This family shows, in a precise way,  
that even small disagreements among dual variables
can cause primal computations to diverge. 
For this reason, we rule out asynchrony in communicating
dual variables. However, we permit all other forms
of asynchrony, and, relative to existing work, this
is the first algorithm to permit arbitrarily asynchronous
computations of dual variables in blocks.

The second contribution of this paper is to establish
convergence rates. 
These rates are shown to depend upon problem parameters,
which lets us calibrate their values 
to improve convergence. Moreover, we show that
convergence can be inexact due to dual asynchrony, and thus
the scalability of parallelization comes at the
expense of a potentially inexact solution. We term
this inexactness the ``asynchrony penalty,'' and we give
an explicit bound on it. 
Simulation results show convergence of this algorithm in practice,
and illustrate concretely that the asynchrony penalty is slight.

The rest of the paper is organized as follows. 
Section~\ref{sec:background} provides the necessary background
on convex optimization and formally gives the asynchronous
primal-dual problem statement. Then Section~\ref{sec:algorithm} discusses four possible asynchronous behaviors, provides a counterexample to complete asynchrony, and presents our asynchronous algorithm. Primal and dual convergence rates are developed in Section~\ref{sec:convergence}. Section~\ref{sec:numerical} presents a numerical example with implications for relationships among parameters. Finally, we present our conclusions in Section~\ref{sec:concl}.

\section{Background and Problem Statement} \label{sec:background}
We study the following form of optimization problem. 

\begin{problem} \label{prob:first}
Given~$h : \R^n \to \R$, ~$g : \R^n \to \R^m$, and~$X \subset \R^n$, asynchronously solve
\begin{align}
\textnormal{minimize }        &h(x) \\
\textnormal{subject to } &g(x) \leq 0 \\
                            &x \in X. \tag*{$\lozenge$}
\end{align} 
\end{problem}

We assume the following about
the objective function~$h$.

\begin{assumption} \label{as:h}
$h$ is twice continuously
differentiable and convex. \hfill $\triangle$
\end{assumption}

We make a similar assumption about
the constraints~$g$.

\begin{assumption} \label{as:g}
$g$ satisfies Slater's condition, i.e.,
there exists~$\bar{x} \in X$ such
that~$g\big(\bar{x}\big) < 0$. For all~$j \in \{1, \ldots, m\}$, the function~$g_j$ is twice continuously differentiable
and convex. \hfill $\triangle$
\end{assumption}

Assumptions~\ref{as:h} and~\ref{as:g} permit
a wide range of functions to be used, such as all
convex polynomials of all orders. 
We impose the following assumption on
the constraint set. 

\begin{assumption} \label{as:X}
$X$ is non-empty, compact, and convex. 
It can be decomposed into~$X = X_1 \times \cdots \times X_N$. 
\hfill $\triangle$
\end{assumption}

Assumption~\ref{as:X} permits many sets to be used,
such as box constraints, which often arise 
multi-agent optimization~\cite{notarnicola2016asynchronous}. 

We will solve Problem~\ref{prob:first}  
using a primal-dual approach. This allows the problem to be parallelized across many processors by re-encoding constraints through Karush-Kuhn-Tucker (KKT) multipliers. In particular, because the constraints~$g$ couple
the processors' computations, they can be difficult to enforce in a distributed way. 
However, by introducing KKT multipliers to re-encode
constraints, we can solve an equivalent, higher-dimensional
unconstrained problem. 

An ordinary primal-dual approach would 
find a saddle point of the Lagrangian associated
with Problem~\ref{prob:first},
defined as~$L(x, \mu) = h(x) + \mu^Tg(x)$.
That is, one would solve~$\min_{x} \max_{\mu} L(x, \mu)$.
However,~$L$ is affine in~$\mu$, which implies
that~$L(x, \cdot)$ is concave but not strongly
concave. Strong convexity has
been shown to provide robustness to asynchrony
in minimization problems~\cite{bertsekas89}, and
thus we wish to endow the maximization 
over~$\mu$ with strong concavity. 
We use a Tikhonov regularization~\cite[Chapter 12]{facchinei2007finite} in~$\mu$
to form
\begin{equation} \label{regL}
L_\delta (x, \mu) = h(x) + \mu^Tg(x) - \frac{\delta}{2} \| \mu \| ^2,
\end{equation}
where~$\delta > 0$. 

Instead of regularizing with respect to the primal variable~$x$, we impose the following assumption
in terms of the Hessian~$H(x,\mu) = \nabla_{x}^2 L_{\delta}(x, \mu)$. When convenient, we suppress
the arguments~$x$ and~$\mu$. 

\begin{assumption}[Diagonal Dominance] \label{as:diagonal}
The~$n \times n$ Hessian matrix~$H=\nabla^{2}_{x}L_{\delta}(x,\mu)$ is~$\beta$-diagonally dominant. That is, for all~$i$ from~$1, \ldots, n$, 
\begin{equation}
|H_{ii}|-\beta \geq \sum_{\substack{ j=1 \\ j \neq i}}^n |H_{ij}|. \tag*{$\triangle$}
\end{equation}
\end{assumption}

If this assumption does not hold, the Lagrangian can be regularized with respect to the primal variable as well, leading to~$H$'s diagonal dominance. Some problems satisfy
this assumption without regularizing~\cite{greene06}, and,
for such problems, we proceed without regularizing
to avoid regularization error. 

Under Assumptions~\ref{as:h}-\ref{as:diagonal},
Problem~\ref{prob:first}
is equivalent to the following
saddle point problem.

\begin{problem} \label{prob:second}
Let Assumptions~\ref{as:h}-\ref{as:diagonal} hold and fix~$\delta > 0$. 
For~$L_{\delta}$ defined in Equation~\eqref{regL}, asynchronously compute
\begin{equation} 
\big(\xhat, \mhatd ) := \argmin_{x \in X} \argmax_{\mu \in \R^m_{+}}  L_\delta (x, \mu). \tag*{$\lozenge$}
\end{equation}
\end{problem}

It is in this form that we will solve
Problem~\ref{prob:first}, and we present our algorithm
for doing so in the next section. 

\section{Asynchronous Primal-Dual Algorithm} \label{sec:algorithm}
One challenge of Problem~\ref{prob:second} is that~$\mhatd$ is maximized over the unbounded domain~$\R^m_{+}$, which
is the non-negative orthant of~$\R^m$. 
Because this domain is unbounded, gradients and other terms are unbounded, which makes convergence analysis challenging as dual iterates may not be within a finite distance of the optimum. 
To remedy this problem, we next compute a non-empty, compact, and
convex set~$M$ that contains~$\mhatd$.

\begin{lemma} \label{lem:mubound}
Let Assumptions~\ref{as:h}-\ref{as:diagonal} hold,
let~$\bar{x}$ be a Slater point of~$g$, and set~$h^* := \min_{x \in X} h(x)$. Then
\begin{equation*}
\mhatd \in M := \Bigg\lbrace	\mu \in \R^m_+ : \| \mu \|_1 \leq \frac{h(\bar{x}) - h^*}{\min\limits_{1 \leq j \leq m}  -g_j(\bar{x})} \Bigg\rbrace.
\end{equation*}
\end{lemma}
\emph{Proof:} 
Follows Section~II-C in \cite{hale15}. \hfill $\blacksquare$

Here,~$h^*$ denotes the optimal unconstrained objective
function value, though any lower-bound for this value will
also provide a valid~$M$. In particular,~$h$ is often non-negative
and one can substitute~$0$ in place of~$h^*$ in such cases. 

Solving Problem~$2$ asynchronously requires choosing an update law that 
we expect to be robust to asynchrony and simple to implement in
a distributed fashion. In this context, first-order
gradient-based methods offer both some degree of inherent
robustness, as well as computations that are simpler than
other available methods, such as Newton-type methods. 
We apply a projected gradient method to both the primal and dual variables,
which is shown in Algorithm~\ref{alg:uzawa}, and is based on the 
seminal Uzawa iteration~\cite{arrow58}. 

\begin{algorithm}[H]
\caption{}
Let~$x$(0) and~$\mu$(0) be given. For values~$k = 0,1,...$, execute
\begin{align}
x(k+1) &= \Pi_X [x(k) - \gamma (\nabla_x L_\delta (x(k),\mu(k)))] \label{eq:uzawax}\\ 
\mu(k+1) &= \Pi_M [\mu(k) + \gamma (\nabla_\mu L_\delta (x(k),\mu(k)))] \label{eq:uzawamu}
\end{align}
where~$\gamma$ is a step-size,~$\Pi_X$ is the Euclidean projection onto~$X$, and~$\Pi_M$ is the Euclidean projection onto~$M$.
\label{alg:uzawa}
\end{algorithm}

\subsection{Overview of Approach}
We are interested in distributing Algorithm \ref{alg:uzawa} among a number of processors while allowing agents to 
generate and share information
as asynchronously as possible. We consider~$N$ agents indexed over~$a \in [N]:=\lbrace 1, \ldots , N\rbrace$. We partition the set~$[N]$ into~$[N_{p}]$ and~$[N_{d}]$
(where~$N_p + N_d = N$). The set~$[N_{p}]$ contains indices of agents that update primal variables (contained in~$x$), while~$[N_{d}]$ contains indices of agents that update dual variables (contained in~$\mu$).
Using a primal-dual approach, there are four behaviors
that could be asynchronous: (i) computations of updates to primal variables, (ii) communications
to share updated values of primal variables, (iii) computations of
updates to dual variables, and (iv) communications to share updated
values of dual variables.

\paragraph*{(i) Computations of Updates to Primal Variables} 
When parallelizing Equation~\eqref{eq:uzawax} across the~$N_p$ primal agents, 
we index all primal agents' computations
using the same iteration counter,~$k \in \N$. 
However, they may perform updates at different times. 
The subset of times at which primal agent~$i \in [N_p]$ computes an update is denoted by~$K^i \subset \N$. 
For distinct~$i, j \in [N_p]$,~$K^i$ and~$K^j$ need not have any relationship. 

\paragraph*{(ii) Communications of Updated Primal Variables} Primal variable communications are also totally asynchronous. A primal variable's current value may be sent to other primal and dual agents that need it at each time~$k$. 
We use the notation~$P^{i}_{j}$ to denote the set of 
times\footnote{We assume that there is no delay between sending and receiving messages. There is no loss of generality in our results because we
can make~$P^{i}_{j}$ and~$D^{i}_{c}$ the times at which messages are received. However, assuming zero delays simplifies the forthcoming discussion and analysis and is done for the remainder
of the paper.} 
at which primal agent~$i$ sends values of its primal variables to agent~$j$. 
Similarly, we use the notation~$D^{i}_{c}$ to denote the set of times at which primal agent~$i$ sends updated values to dual agent~$c \in [N_d]$.

\paragraph*{(iii) Computations of Updates to Dual Variables} 
Dual agents wait for each primal agent's updated state before computing an update. 
Dual agents may perform updates at different times because they may receive primal updates at different times. In some cases, a dual agent may receive multiple updates from a subset 
of primal agents prior to receiving all required primal updates. In this case, only the most recently received update from a primal agent will be used in the dual agent's computation.
For all~$c \in [N_d]$, dual agent~$c$ keeps an iteration count~$t_c$ to track the number of updates it has completed. 

\paragraph*{(iv) Communications of Updated Dual Variables} 
Previous work~\cite{hale17} has shown that allowing primal agents to disagree arbitrarily about dual variables
can outright preclude convergence. 
This is explained by the following: fix~$\mu^1, \mu^2 \in M$. Then an agent with~$\mu^1$ onboard
is minimizing~$L(\cdot, \mu^1)$, while an agent with~$\mu^2$ onboard is minimizing~$L(\cdot, \mu^2)$. 
If~$\mu^1$ and~$\mu^2$ are arbitrarily far apart, then it is not surprising that
the minima of~$L(\cdot, \mu^1)$ and~$L(\cdot, \mu^2)$ are as well. 
However, one may conjecture that small disagreements in dual variables lead to small
distances between these minima. 
Below, we show that this conjecture is false and that even small disagreements in dual variables can lead to arbitrarily 
large distances between the minima they induce. Even limited asynchrony can lead to small disagreements in dual
variables, and, in light of the above discussion, this can cause primal agents' computations to reach points
that are arbitrarily far apart. 

Therefore, we will develop an algorithm that proceeds
with all agents agreeing on the value of~$\mu$ while still allowing 
dual
computations to be divided among dual agents. This is accomplished by allowing primal agents to work completely asynchronously (updates are computed and sent at different times to different agents) but requiring that dual updates are sent to all primal agents at the same time\footnote{Even if they are not sent and/or received at the same time, we can apply any procedure to synchronize these values and the algorithm will remain the same. We assume synchrony in sending and receiving these values merely
to simplify the forthcoming discussion.}. 
After a dual agent computes an update, it sends its updated dual variable to all primal agents. 
Again, for simplicity, it is assumed that the update is received at the same time it is sent.   When dual agent~$c$ sends the updated dual variable~$\mu^{c}_{c}$ to all primal agents, it also sends its iteration count~$t_{c}$. This allows primal agents to annotate which version of~$\mu$ is used in their updates. Primal agents 
disregard any received primal updates that use an outdated version of~$\mu_c$ (as indicated by~$t_c$). 
This ensures that primal updates are not mixed if they rely on different dual values. 

\subsection{Counterexample to the Asynchronous Dual Case} \label{sec:counterexample}
Below we show 
that behavior (iv) above, communications to share updated values
of dual variables, cannot be asynchronous in general. 
The intuition here is as follows. 
In a primal-dual setup, one can regard each fixed choice of dual vector as
specifying a problem to solve in a parameterized family of minimization problems.
Formally, with~$\mu$ fixed, agents solve 
\begin{equation} \label{eq:mingradl}
\textnormal{minimize}_{x \in X} L_{\delta}(x, \mu) := h(x) + \mu^Tg(x) - \frac{\delta}{2}\|\mu\|^2. 
\end{equation}
For two primal agents with different values
of~$\mu$, denoted~$\mu^1$ and~$\mu^2$, they solve two different
problems: agent~$1$ minimizes~$L_{\delta}(\, \cdot\,, \mu^1)$ while
agent~$2$ minimizes~$L_{\delta}(\,\cdot\,, \mu^2)$. 
With a gradient-based method to minimize over~$x$,
gradients depend linearly upon~$\mu$. 
This may 
lead one to believe that for
\begin{equation}
\hat{x}_1 := \argmin_{x \in X} L_{\delta}(x, \mu^1) \quad \textnormal{ and } \quad
\hat{x}_2 := \argmin_{x \in X} L_{\delta}(x, \mu^2),
\end{equation}
having~$\|\mu^1 - \mu^2\|$ small implies
that~$\|\hat{x}_1 - \hat{x}_2\|$ is also small. However, we show in the following theorem
that this is false. 

\begin{theorem} \label{thm:varyingmu}
Fix any~$\epsilon > 0$ and~$L > \epsilon$. Then, under Assumptions 1-4, there always exists
a problem such that~$\|\mu^1 - \mu^2\| < \epsilon$ and~$\|\hat{x}_1 - \hat{x}_2\| > L$. 
\end{theorem}
\emph{Proof:} See the appendix.\hfill $\blacksquare$

\subsection{Glossary of Notation}
Every agents stores a local copy of~$x$ and~$\mu$ for use in local computations. 
The following notation is used in our formal algorithm statement below. 
\begin{itemize}

\item[$D^{i}_{c}$] The times at which messages are sent by primal agent~$i$ and received by dual agent~$c$. 
\item[$g_c(x)$] The~$c^{th}$ entry of the constraint function,~$g$, evaluated at~$x$.
\item[$k$] The iteration count used by all primal agents.
\item[$K^{i}$] The set of times at which primal agent~$i$ performs updates.
\item[$\mathcal{N}_i$] Essential neighborhood of agent~$i$. Agent~$j$ is an essential neighbor of agent~$i$ if $\nabla_{x_{j}}L_{\delta}$ depends upon~$x_{i}$. Then agent~$i$ communicates with agent~$j$ to ensure
it has the information necessary to compute gradients.
\item[${[N_{d}]}$] Set containing the indices of all dual agents.
\item[${[N_{p}]}$] Set containing the indices of all primal agents.
\item[$P^{i}_{j}$] The times at which messages are sent by primal agent~$i$ and received by primal agent~$j$. 
\item[$\tau^{i}_{j}(k)$] Time at which primal agent~$j$ computed the update that it sent agent~$i$ at time~$k$ ($i$ can be primal or dual). Note that~$\tau^{i}_{i}(k)=k$ for all~$i \in [N_p]$.
\item[$t$] The vector of dual agent iteration counts. The~$c^{th}$ entry,~$t_c$, is the iteration count for dual agent~$c$'s updates. 
\item[$t_{c}$] The iteration count for dual agent~$c$'s updates. This is sent along with~$\mu^{c}_{c}$ to all agents.
\item[$v^{c}_{i}$] The iteration count used by dual agent~$c$ for updates received from primal agent~$i$ between dual updates.
\item[$x^{i}_{j}$] Agent~$i$'s value for the primal variable~$j$, which is updated/sent by primal agent~$j$. If agent~$i$ is primal, it is indexed by both~$k$ and~$t$; if agent~$i$ is dual it is indexed by~$t$.
\item[$x^{i}_{i} (k;t)$] Agent~$i$'s value for its primal variable~$i$ at primal time~$k$, calculated with dual update~$t$.
\item[$x^*(t)$] The fixed point of~$f=\Pi_{X}\left[x - \gamma\nabla_{x}L_{\delta}(x,\mu(t))\right]$ with respect to a fixed~$\mu(t)$.
\item[$x^c_t$] Abbreviation for $x^c(t)$, which is dual agent~$c$'s copy of the primal vector at time~$t$.
\item[$\xhat$] The primal component of the saddle point of~$L_{\delta}$. Part of the optimal solution pair~$(\xhat , \mhatd)$.
\item[$\xhatt$] Given~$\mu(t)$,~$\xhatt = \argmin_{x \in X} L_{\delta}(x, \mu (t))$.
\item[$\mu^{c}_{d}$] Agent~$c$'s copy of dual variable~$d$, which is updated/sent by dual agent~$d$. Agent~$c$ may be primal or dual.
\item[$\mhatd$] The dual component of the saddle point of~$L_{\delta}$,~$(\xhat , \mhatd)$.
\item[$\mhat$] The~$c^{th}$ entry of~$\mhatd$.
\item[$M_c$] The set~$\{\nu \in \mathbb{R}_{+} : \nu\leq\frac{h(\bar{x}) - h^*}{\min_{j} -g_j(\bar{x})}\}$. This uses the upper bound in Lemma~\ref{lem:mubound} to project
individual components of~$\mu$. 
\end{itemize}

\subsection{Statement of Algorithm}
Having defined our notation, we impose the following assumption on agents' communications and computations.
\begin{assumption} \label{as:times}
For all~$i \in [N_p]$, the set~$K^i$ is infinite. If~$\{k_n\}_{n \in \N}$ is an increasing set
of times in~$K^i$, then~$\lim_{n \to \infty} \tau^j_i(k_n) = \infty$ for all~$j \in [N_p]$
and~$\lim_{n \to \infty} \tau^c_i(k_n) = \infty$ for all~$c \in [N_d]$. \hfill $\triangle$
\end{assumption}
This simply ensures that no agent ever stop computing or communicating, though
delays can be arbitrarily large. 

We now define the asynchronous primal-dual algorithm. 

\begin{algorithm}[H]
\caption{}
Step 0: Initialize all primal and dual agents with~$x(0)\in X$ and~$\mu(0)\in M$. Set~$t=0$ and~$k=0$. \\
Step 1: For all~$i\in[N_{p}]$ and all~$j\in \mathcal{N}_{i}$, if~$k \in P^{i}_{j}$, then agent~$j$ sends~$x^{j}_{j}(k;t)$ to agent~$i$.\\
Step 2: For all~$i\in[N_{p}]$ and all~$c\in[N_{d}]$, if agent~$i$ receives a dual variable update from agent~$c$, it uses the accompanying~$t_{c}$ to update the vector~$t$ and performs the update
\begin{equation*}
\mu^{i}_{c}(t) = \mu^{c}_{c}(t_{c}).
\end{equation*}
Step 3: For all~$i\in[N_{p}]$ and all~$j\in \mathcal{N}_{i}$, execute
\begin{align*}
x^{i}_{i}(k\!+\!1;t) \!&=\! \begin{cases} 
                              \Pi_{X_{i}} [x^{i}_{i}(k;t) \!-\! \gamma\nabla_{x_{i}}L_{\delta}(x^{i}(k;t),\mu(t))] & k \!\in\! K^{i}\\
                              x^{i}_{i}(k;t) & k \!\notin\! K^{i} 
                              \end{cases} \\
x^{i}_{j}(k+1;t) &= \begin{cases} 
                              x^{j}_{j}(\tau^{i}_{j}(k+1);t) & i$ receives~$j$'s state at~$k+1 \\
                              x^{i}_{j}(k;t) & $otherwise$
                              \end{cases}.                    
\end{align*}
Step 4: If~$k+1 \in D^{i}_{c}$, agent~$i$ sends~$x^{i}_{i}(k+1;t)$ to dual agent~$c$. Set~$k:=k+1$.\\
Step 5: For~$c \in [N_{d}]$ and~$i\in[N_{p}]$, if dual agent~$c$ receives an update from primal agent~$i$ computed with dual update~$t$, it sets
\begin{align*}
x^{c}_{i}(t_{c}) &= \begin{cases} 
                              x^{i}_{i}(\tau^{c}_{i}(k),t) & x^{i}_{i}$ received prior to update~$t_{c}+1 \\
                              x^{c}_{i}(t_{c}-1) & $otherwise$
                               \end{cases}.
\end{align*}
Step 6: For~$c \in [N_{d}]$, if agent~$c$ has received an update from every primal agent for the latest dual iteration~$t$, it executes
\begin{equation*}
\mu^{c}_{c}(t_{c}+1) = \Pi_{M_{c}}[\mu^{c}_{c}(t_{c}) + \rho\frac{\partial L}{\partial \mu_{c}}(x^{c}(t_{c}),\mu^{c}(t_{c}))] .
\end{equation*}
Step 7: If dual agent~$c$ updated in Step 6, then it sends~$\mu^{c}_{c}(t_c + 1)$ to all primal agents. Set~$t_{c}:=t_{c}+1$.\\
Step 8: Return to Step 1.
\label{alg:2}
\end{algorithm}

\section{Convergence} \label{sec:convergence}

To define an overall convergence rate to the optimal solution~$(\xhat , \mhatd)$, we first fix the dual variable~$\mu$ and find the primal convergence rate. We then find the overall dual convergence rate to~$\mhatd$ 
by showing that dual variables converge to~$\mhatd$ over time, which lets us show that primal variables converge to~$\xhat$. 

\subsection{Primal Convergence with Fixed Dual Variable}
Given a fixed~$\mu(t)$, projected gradient descent for minimizing~$L_{\delta}(\cdot,\mu(t))$ may be written as
\begin{equation} \label{eq:grad}
f(x) = \Pi_{X}\left[x - \gamma\nabla_{x}L_{\delta}(x,\mu(t))\right], 
\end{equation}
where~$\gamma > 0$. 
Leveraging some existing theoretical tools in the study of
optimization algorithms~\cite{Bertsekas1983, Bertsekas1991}, we can
study~$f$ in a way that elucidates its behavior under asynchrony.

According to \cite{Bertsekas1991}, the assumption of diagonal dominance guarantees that~$f$ has the contraction property 
\begin{equation}
{\Vert f(x)-x^*(t) \Vert \leq \alpha \Vert x - x^*(t) \Vert} \textnormal{ for all } x \in X, 
\end{equation}
where~$\alpha \in [0,1)$ and~$x^*(t)$ is a fixed point
of~$f$, which depends on the choice of fixed~$\mu(t)$.  
However, 
the value of~$\alpha$ is not specified in~\cite{Bertsekas1991}, 
and it is precisely that value that governs the rate of convergence
to a solution. We therefore compute~$\alpha$ explicitly below.
First, we bound the step-size~$\gamma$. 

\begin{definition} \label{def:gamma}
Define the primal step-size~$\gamma > 0$ to satisfy
\begin{equation}
\gamma < \frac{1}{\max\limits_i\max\limits_{x \in X} \max\limits_{\mu \in M} \sum_{j=1}^n |H_{ij}(x, \mu)|} .
\end{equation}
\end{definition}
Because~$x$ and~$\mu$ both take values in compact sets, each entry~$H_{ij}$ is bounded
above and below, and thus the upper bound on~$\gamma$ is positive. 

Following the method in \cite{Bertsekas1983}, two~$n \times n$ matrices~$G$ and~$F$ must also be defined.
\begin{definition} Define the~$n \times n$ matrices~$G$ and~$F$ as
\begin{equation}
   G=
    \begin{bmatrix}
    |H_{11}| & -|H_{12}| & \ldots & -|H_{1n}| \\
    \vdots & \vdots & \ddots & \vdots  \\
    -|H_{n1}| & -|H_{n2}| & \ldots & |H_{nn}| \\
    \end{bmatrix} \,\, \textnormal{and} \,\, F = I - \gamma G,
\end{equation}
where~$I$ is the~$n \times n$ identity matrix.
\end{definition}

We now state the following lemma that relies on meeting the conditions listed in~\cite{Bertsekas1983}.

\begin{lemma} \label{83lemma}
Let~$f$,~$G$, and~$F$ be as above and let Assumptions~\ref{as:h}-\ref{as:diagonal} hold. 
Then,
$|f(x) - f(y)| \leq F|x-y|$,
 for all~$x,y \in \R^n$, 
where~$|v|$ denotes the element-wise absolute value
of the vector~$v \in \mathbb{R}^n$ and the inequality holds component-wise. 
\end{lemma}
\emph{Proof:}
Three conditions must be satisfied in \cite{Bertsekas1983}: (i)~$\gamma$ is sufficiently small, (ii)~$G$ is positive definite, and (iii)~$F$ is positive definite.

\paragraph*{(i)~$\gamma$ is sufficiently small} 
Results in \cite{Bertsekas1991} require
$\gamma \sum_{j=1}^n |H_{ij}| < 1$  for all $i \in \lbrace 1,\ldots, n\rbrace$,
which here follows immediately from Definition~\ref{def:gamma}.

\paragraph*{(ii)~$G$ is positive definite}
By definition,~$G$ has only positive diagonal entries. By~$H$'s diagonal dominance we have the following inequality for all~$i \in \lbrace 1,\ldots, n\rbrace$:
\begin{equation}
|G_{ii}| = |H_{ii}| \geq \sum_{\substack{ j=1 \\ j \neq i}}^n |H_{ij}| + \beta > \sum_{\substack{ j=1 \\ j \neq i}}^n |H_{ij}| = \sum_{\substack{ j=1 \\ j \neq i}}^n |G_{ij}|.
\end{equation}
Because~$G$ has positive diagonal entries, is symmetric, and is strictly diagonally dominant,~$G$ is positive definite by Gershgorin's Circle Theorem. 

\paragraph*{(iii)~$F$ is positive definite}
Definition~\ref{def:gamma} ensures the diagonal entries of~$F$ are always positive. And $F$ is diagonally dominant if, for all~$i \in \{1, \ldots, n\}$,
\begin{equation}
|F_{ii}| = 1-\gamma |H_{ii}| > \gamma \sum_{\substack{ j=1 \\ j \neq i}}^n |H_{ij}| = \sum_{\substack{ j=1 \\ j \neq i}}^n |F_{ij}| .
\end{equation}
This can be rewritten as~$\gamma \sum_{j=1}^n |H_{ij}|  < 1$, 
which was satisfied under Condition~(i). 
Because~$F$ has positive diagonal entries, is symmetric, and is strictly diagonally dominant,~$F$ is positive definite by Gershgorin's Circle Theorem.\hfill $\blacksquare$

To establish convergence in~$x$, we will show that~${\Vert f(x)-f(x^*(t)) \Vert \leq \alpha \Vert x - x^*(t) \Vert}$ for all~$x \in X$, where~$\alpha \in [0,1)$. 
Toward doing so, we define the following.
\begin{definition} \label{def:vmax}
Let $v \in \R^n$. Then~$\Vert v\Vert _{max} = \max\limits_i |v_i|$. 
\end{definition}
In this work, we consider a scalar maximum norm because we consider
agents that update scalar blocks, though updating
non-scalar blocks is readily accommodated by considering a block-maximum norm~\cite{bertsekas89}. 

We next show that the gradient update law~$f$ in Equation~\eqref{eq:grad} converges with asynchronous, distributed computations. Furthermore, we quantify convergence in terms of~$\gamma$ and~$\beta$. 

\begin{lemma} \label{lem:primal}
 Let~$f$,~$H$,~$G$,~$F$,~$\gamma$, and~$\Vert v\Vert _{max}$ be as defined above. Let Assumptions~\ref{as:h}-\ref{as:diagonal} hold and fix~$\mu(t) \in M$. 
Then for a fixed point~$x^*(t)$ of~$f$ and for all~$x \in \R^n$, 
\begin{equation}
\Vert f(x)-f(x^*(t)) \Vert _{max} \leq q_p\Vert x-x^*(t)\Vert _{max},
\end{equation}
where~$q_p=(1-\gamma \beta) \in [0,1)$. 
\end{lemma}

\emph{Proof:}
For notational simplicity we write~$x^*(t)$ simply as~$x^*$. 
Assumption~\ref{as:diagonal} and the definition of~$F$ give
\begin{equation}
\sum_{\substack{ j=1}}^n F_{ij} = 1- \gamma\Big(|H_{ii}| - \sum_{\substack{ j=1 \\ j \neq i}}^n |H_{ij}|\Big) \leq 1- \gamma\beta.
\end{equation}

This result, Definition~\ref{def:vmax}, and Lemma~\ref{83lemma} give
\begin{align*}
&\Vert f(x)-f(x^*) \Vert _{max} = \max\limits_i |f_i(x) - f_i(x^*)| \\
&\quad\leq \max\limits_i \sum_{\substack{ j=1}}^n F_{ij} |x_j - x_j^*| \leq \max\limits_l |x_l - x_l^*| \max\limits_i \sum_{\substack{ j=1}}^n F_{ij} \\
&\quad\leq \max\limits_l |x_l - x_l^*| (1- \gamma \beta) = (1-\gamma \beta)\Vert x-x^*\Vert _{max}. 
\end{align*}
All that remains is to show~$(1-\gamma \beta) \in [0,1)$. By Definition~\ref{def:gamma} and the inequality~$|H_{ii}|\geq \beta$, for all~$x \in X$ and~$\mu(t) \in M$,
\begin{equation}
\gamma \beta < \frac{\beta}{\max\limits_i \sum_{j=1}^n |H_{ij}|} \leq \frac{\beta}{\max\limits_i |H_{ii}|} < \frac{\beta}{\beta} = 1. \tag*{$\blacksquare$} 
\end{equation}

The primal-only convergence rate can be computed by
leveraging results in~\cite{hale17} in terms of the number
of operations the primal agents have completed (counted in the appropriate
sequence). Namely, we count operations as follows. For a given dual variable with iteration vector~$t$, we set~$\textnormal{ops}(k,t) = 0$. Then, after
all primal agents have computed an update to their decision variable and
sent it to and had it received by all other primal agents
that need it, say by time~$k'$, we increment~$\textnormal{ops}$
to~$\textnormal{ops}(k',t) = 1$. After~$\textnormal{ops}(k',t) = 1$,
we then wait until all primal agents have subsequently computed
a new update (still using the same dual variable indexed with~$t$) and it has been received by all other primal agents
that need it. If this occurs at time~$k''$, then
we set~$\textnormal{ops}(k'',t) = 2$, and then this process
continues. If at some time~$k'''$, primal agents receive an updated~$\mu$ (whether just a single dual agent sent an update or multiple agents send updates) with an iteration vector of~$t'$, then the cycle count would begin again with~$\textnormal{ops}(k''',t') = 0$.


\begin{theorem} \label{thm:primalconv}
Let Assumptions~\ref{as:h}-\ref{as:diagonal} hold. 
For~$\mu(t)$ fixed and the agents asynchronously executing
the gradient update law~$f$,
\begin{equation}
\|x^i(k) - x^*(t)\|_{max} \leq q_p^{\textnormal{ops}(k,t)} \max\limits_{j} \|x^j(k) - x^*(t)\|_{max}, 
\end{equation}
where~$x^*(t)$ is the fixed point of~$f$ with~$\mu(t)$ fixed. 
\end{theorem}
\emph{Proof:} 
We write~$x^*$ in place of~$x^*(t)$ for simplicity. 
From Lemma~\ref{lem:primal} we see that~$f$ is
a~$q_p$-contraction mapping with respect to the norm~$\|\cdot\|_{max}$. 
From Section 6.3 in~\cite{Bertsekas1991}, this
contraction property implies that there
exist sets of the form
\begin{equation}
X(k) = \{x \in \R^n \mid \|x - x^*\|_{max} \leq q_p^k\|x(0) - x^*\|_{max}\}
\end{equation}
that satisfy the following criteria from~\cite{hale17}:
\begin{enumerate}[i.]
\item $\cdots \subset X(k+1) \subset X(k) \subset \cdots \subset X$
\item $\lim_{k \to \infty} X(k) = \{x^*\}$
\item For all~$i$, there are sets~$X_i(k) \subset X_i$
satisfying
\begin{equation}
X(k) = X_1(k) \times \cdots \times X_N(k)
\end{equation}
\item For all~$y \in X(k)$ and all~$i \in [N_p]$, 
$f_i(y) \in X_i(k+1)$, where
$f_i(y) = \Pi_{X_i}\left[y_i - \gamma \nabla_{x_i} L_{\delta}(y)\right]$. 
\end{enumerate}

We will use these properties to compute the desired
convergence rate. Suppose all agents have
a fixed~$\mu(t)$ onboard. Upon receipt of this~$\mu(t)$,
agent~$i$ has~$x^i(k; t) \in X(0)$ by definition.
Suppose at time~$\ell_i$ that agent~$i$ computes
a state update. Then~$x^i_i(\ell_i+1; t) \in X_i(1)$.
For~$m = \max_{i \in [N_p]} \ell_i + 1$, we find
that~$x^i_i(m; t) \in X_i(1)$ for all~$i$.
Next, suppose that, after all updates have been computed, 
these updated values are sent to and received by all agents
that need them, say at time~$m'$. Then, for any~$i \in [N_p]$, agent~$i$
has~$x^i_j(m'; t) \in X_j(1)$ for all~$j \in [N_p]$.
In particular,~$x^i(m'; t) \in X(1)$, and this is satisfied
precisely when a single cycle has occurred. Iterating this argument
for subsequent cycles completes the proof. 
\hfill $\blacksquare$

\subsection{Dual convergence}

We next derive a componentwise
convergence rate for the dual variable. 

\begin{theorem} \label{thm:dualsingle}
Let Assumptions~\ref{as:h}-\ref{as:diagonal} hold. 
Fix~$\delta >0$ and let 
$\rho \in \left(\frac{3 - \sqrt{3}}{3 \delta} , \frac{3 + \sqrt{3}}{3 \delta}\right)$.
Then 
\begin{align*}
|\mcc (&t_c+1) - \mhat |^2 \leq q_d |\mcc (t_c) - \mhat |^2 + 2 \rho^2 M_{g_c}^2 D_x^2  + 2 \rho^2 M_{g_c}^2 q_p^{2 \textnormal{ops}(k_c, t)}L_x^2 + 2\rho ^2 M_{g_c}^2 D_x q_p^{\textnormal{ops}(k_c, t)}L_x, 
\end{align*}
where ~${q_d:=3(1- \rho \delta)^2 \in [0,1)}$,~${M_{g_c} := \max\limits_{x \in X} \| \nabla g_c(x) \|_{max}}$, ${D_x := \max\limits_{x,y \in X} \| x-y \|_{max}}$, ${L_x=\max\limits_{j} \|x^j(0) - x^*(t)\|_{max}}$,
$t$ is the dual iteration count vector onboard the primal agents when sending states to dual agent~$c$,
and $k_c$ is the time at which the first primal agent sent a state to dual agent~$c$ with~$\mu(t)$ onboard. 
\end{theorem}

\emph{Proof:} 
Let~$x^c(t)$ be denoted by~$x^c_t$, and define~$\xhatt = \argmin_{x \in X} L_{\delta}(x, \mu(t))$ and~$\xhat = \argmin_{x \in X} L_{\delta}(x, \mhatd)$.

Using the non-expansiveness of~$\Pi_M$, we find
\begin{align*}
|\mcc (&t_c+1) - \mhat |^2 = |\Pi _M [\mcc (t_c) + \rho (g_c(x^c_t) 
    -\delta \mcc (t_c))] - \Pi _M [\mhat + \rho (g_c(\xhat)-\delta \mhat )]|^2 \\
&\leq (1- \rho \delta)^2|\mcc (t_c) - \mhat |^2 + \rho ^2|g_c(x^c_t) - g_c(\xhat)|^2 -2 \rho (1-\rho \delta)(\mcc (t_c) - \mhat )(g_c(\xhat)-g_c(x^c_t)).
\end{align*}

Adding~$g_c(\xhatt) - g_c(\xhatt)$ in the last set of parentheses gives
\begin{align}
|\mcc (t_c+1) - \mhat |^2 &\leq (1- \rho \delta)^2|\mcc (t_c) - \mhat |^2 \nonumber + \rho ^2|g_c(x^c_t) - g_c(\xhat)|^2 \nonumber \\
&\quad -2 \rho (1-\rho \delta)(\mcc (t_c) - \mhat )(g_c(\xhat)-g_c(\xhatt)) \nonumber \\
& \quad-2 \rho (1-\rho \delta)(\mcc (t_c) - \mhat )(g_c(\xhatt)-g_c(x^c_t)). 
\label{eq:addsub}
\end{align}

Using
${0 \leq |(1-\rho \delta) (\mcc (t_c) - \mhat ) + \rho (g_c(\xhatt)-g_c(x^c_t))|^2}$,
we expand and rearrange to give
\begin{multline} \label{eq:bigineq1}
{-2}\rho (1-\rho \delta) (\mcc (t_c) - \mhat) (g_c(\xhatt)-g_c(x^c_t))
\leq (1-\rho \delta)^2 |\mcc (t_c) - \mhat |^2  + \rho ^2 |g_c(\xhatt)-g_c(x^c_t)|^2.
\end{multline}
Similarly, we can derive
\begin{multline} \label{eq:bigineq2}
{-2}\rho (1-\rho \delta) (\mcc (t_c) - \mhat) (g_c(\xhat)-g_c(\xhatt))
\leq (1- \rho \delta) ^2 |\mcc (t_c) - \mhat|^2  + \rho ^2 |g_c(\xhat)-g_c(\xhatt)|^2.
\end{multline}

Using Equations~\eqref{eq:bigineq1} and~\eqref{eq:bigineq2} in Equation~\eqref{eq:addsub} gives
\begin{align}
|\mcc (t_c\!\!+\!\!1) \!-\! \mhat |^2 &\!\!\leq\! 3(1 \!-\! \rho \delta)^2 |\mcc (t_c) \!-\! \mhat |^2 +\!\!\rho ^2 |g_c(\xhat)\!-\!g_c(\xhatt)|^2 \\
&\quad \!\!+\! \rho ^2 |g_c(\xhatt)\!-\!g_c(x^c_t)|^2 +\! \rho ^2|g_c(x^c_t) \!-\! g_c(\xhat)|^2.
\label{eq:grads}
\end{align}
For the~$\rho ^2|g_c(x^c_t) - g_c(\xhat)|^2$ term, we can write
\begin{align*}
|g_c(x^c_t) - g_c(\xhat)|^2 &= |g_c(x^c_t) - g_c(\xhatt) +g_c(\xhatt) - g_c(\xhat)|^2 \\
&\leq |g_c(x^c_t) - g_c(\xhatt)|^2 + |g_c(\xhatt) - g_c(\xhat)|^2 + 2|g_c(x^c_t) - g_c(\xhatt)| |g_c(\xhatt) - g_c(\xhat)|,
\end{align*}
where substituting this into Equation~\eqref{eq:grads} and grouping gives
\begin{align}
|\mcc (t_c+1) - \mhat |^2 &\leq 3(1- \rho \delta)^2 |\mcc (t_c) - \mhat |^2  \nonumber + 2 \rho^2 |g_c(\xhatt) - g_c(\xhat)|^2 \nonumber \\
&\quad + 2 \rho^2 |g_c(x^c_t) - g_c(\xhatt)|^2 + 2\rho ^2|g_c(x^c_t) - g_c(\xhatt)| |g_c(\xhatt) - g_c(\xhat)|.
\end{align}
Using the Lipschitz property of~$g_c$ and the definition of~$D_x$,
\begin{align}
|\mcc (t_c\!+\!1) \!-\! \mhat |^2 &\leq 3(1\!-\! \rho \delta)^2 |\mcc (t_c) \!-\! \mhat |^2 \!+\! 2 \rho^2 M_{g_c}^2 D_x^2 \nonumber \\
&\quad +2 \rho^2 M_{g_c}^2 \|x^c_t - \xhatt\|_{max}^2  + 2\rho ^2 M_{g_c}^2 D_x\|x^c_t - \xhatt\|_{max}. \label{eq:xgrads}
\end{align}
Previously, we established primal convergence for a fixed~$\mu(t)$. This allows us to write
$\|x^c_t - \xhatt\|_{max} \leq q_p^{\textnormal{ops}(k_c,t)}L_x$.
Substituting this into Equation~\eqref{eq:xgrads} completes the proof. \hfill $\blacksquare$

\begin{remark}
The term~$2\rho^2M_{g_c}^2D_x^2$ in Theorem~\ref{thm:dualsingle} is termed the ``asynchrony penalty,'' because
it is an offset from reaching a solution. It is due to asynchronously computing dual
variables and is absent when dual updates are centralized~\cite{koshal2011multiuser,hale17}. 
Further bounding it will be the subject of future research. 
\end{remark}

We next present a simplified dual convergence rate. 

\begin{theorem} \label{thm:geometric}
Let all conditions of Theorem~\ref{thm:dualsingle} hold. 
In Algorithm \ref{alg:2}, convergence for dual agent~$c$ obeys
\begin{equation*}
    | \mcc (t_c+1) - \mhat |^2 \leq q_d^{t_c+1} | \mcc (0) - \mhat |^2 + \frac{1-q_d^{t_c+2}}{1-q_d} p ,
\end{equation*}
where~$p = \max\limits_t 2 \rho^2 M_{g_c}^2 D_x^2+2 \rho^2 M_{g_c}^2 q_p^{2 \textnormal{ops}(k_c,t)}L_x^2 +2\rho ^2 M_{g_c}^2 D_x q_p^{\textnormal{ops}(k_c,t)}L_x.$
\end{theorem}

\emph{Proof:} First, define 
\begin{equation*}
p_i := 2 \rho^2 M_{g_c}^2 D_x^2 + 2 \rho^2 M_{g_c}^2 q_p^{2 \textnormal{ops}(k_c,i)}L_x^2 + 2\rho ^2 M_{g_c}^2 D_x q_p^{\textnormal{ops}(k_c,i)}.
\end{equation*}
Then, recursively applying Theorem~\ref{thm:dualsingle} gives
\begin{equation}
    | \mcc (t_c+1) \!-\! \mhat |^2 \!\leq\! q_d^{t_c+1} | \mcc (0) - \mhat |^2 + \sum_{j=0}^{t_c + 1}q_d^{t_c+1-j} p_j. 
    \label{eq:sum}
\end{equation}
By definition,~$q_d \in [0,1)$, so, using~$p_i \leq p$, summing the geometric series in Equation~\eqref{eq:sum} 
completes the proof. \hfill $\blacksquare$

Theorems~\ref{thm:dualsingle} and~\ref{thm:geometric} can be used to give an overall primal convergence rate
for~$x^i(k; t)$ converging to~$\hat{x}_{\delta}$.

\begin{theorem}
Let Assumptions~\ref{as:h}-\ref{as:diagonal} hold. Then
there exist constants~$K_1, K_2 > 0$ such that, for all~$t$,
\begin{equation}
\|x^i(k; t) -  \hat{x}_{\delta}\|_2 \leq K_1q_p^{\textnormal{ops}(k,t)} + K_2\|\mu(t) - \hat{\mu}_{\delta}\|_2.
\end{equation}
\end{theorem}
\emph{Proof:} Plug Theorem~\ref{thm:dualsingle} into
Theorem~2 of~\cite{hale17}. \hfill $\blacksquare$

\section{Numerical Example} \label{sec:numerical}
We consider an example with~$n = 10$ primal agents (each updating a scalar variable) whose objective function is
\begin{equation}
h(x) = \sum_{i=1}^{n} x_i^4 + \frac{1}{20}\sum_{i=1}^{n} \sum_{\substack{j=1 \\ j \neq i}}^{n} (x_i - x_j)^2.
\end{equation}

Set~$\delta = 0.001$. For~$b = (-2, 4, -10, 5, 1, 8)^T$ and
\begin{equation*}
A = \left(\begin{array}{cccccccccc}
-1&0 & -3 & 0 & 0 & 4 & 0 & 0 & 10 & 0 \\
0 & 1 & 5 & 1 & 1 & 0 & 0 & 2 & 0 & 5 \\
0 & 0 & 1 & 1 & -5 & 1 & 4 & 0 & 0 & 0 \\
0 & 0 & -2 & 0 & 0 & 8 & 1 & 1 & -3 & 1 \\
0 & 0 & 0 & 0 & -3 & 0 & 1 & 1 & 1 & 0 \\
0 & 4 & 0 & 0 & 0 & 0 & 0 & 2 & 1 & -4
\end{array}\right),
\end{equation*}
there are~$m=6$ dual agents each responsible for a scalar dual variable that encodes a constraint in~$g(x) = Ax - b$. 
We also require that each~$x_i \in X_i = [1, 10].$
Then
\begin{equation*}
L_{\delta}(\!x,\!\mu)
=  \sum_{i=1}^{n} x_i^4 + \frac{1}{20}\sum_{i=1}^{n} \sum_{\substack{j=1 \\ j \neq i}}^{n} (x_i\!-\!x_j)^2
+ \mu^T\Big(\!Ax\!-\!b\Big) - \frac{\delta}{2}\|\mu\|^2. 
\end{equation*}
We find that~$H$ is~$\beta$-diagonally dominant with~$\beta = 12$. 

We use this simulation example to explore fundamental relationships between the diagonal dominance parameter~$\beta$ and the communication rate between agents. Here, we use
random communications and vary the probability that agents~$i$ and~$j$ communicate at a particular timestep. 
We begin by varying~$\beta$  over~$\beta \in 12 \cdot \lbrace.9,1,10,100\rbrace$, which we do by scaling~$h$ while holding other terms constant.  
 Figure~1 plots the iteration number~$k$ versus the value~$\frac{\Vert x^i_i(k; t) - \xhat \Vert}{\Vert \xhat \Vert}$ (the relative error) for these values of~$\beta$ and a communication rate of~$1$ (the agents communicate every update with one another).

\begin{figure}[t!]
\centering
\includegraphics[width=8.8cm]{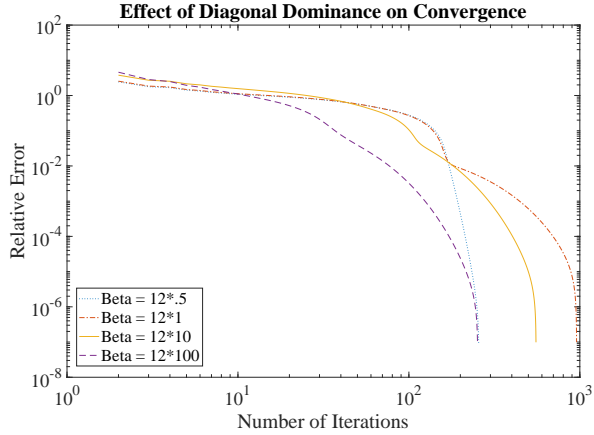}
\caption{Effect of diagonal dominance on convergence. Larger values of~$\beta$ tend to produce faster convergence. However, sufficiently small values for~$\beta$ may also produce the same result by promoting constraint satisfaction earlier.}
\label{fig}
\end{figure}
\begin{figure}[t!]
\centering
\includegraphics[width=8.8cm]{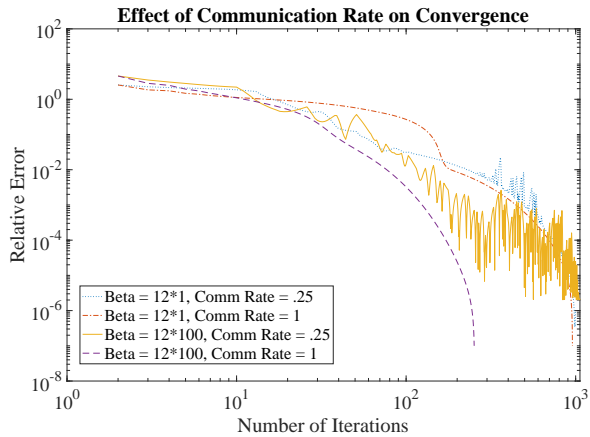}
\caption{Effect of communication rate on convergence. Less frequent communication leads to slower convergence, larger errors, and oscillation in the optimum's proximity. The benefits of strong diagonal dominance can also be outweighed by poor communication.}
\label{fig2}
\end{figure}

As predicted by Theorem~\ref{thm:primalconv}, a larger~$\beta$ correlates with faster convergence in general. Figure~1 also reveals, however, that if~$\beta$ is sufficiently small, then the convergence rate is 
eventually similar to that for a large~$\beta$. 
This is explained by how~$\beta$ weights~$h$ versus~$g$: if~$\beta$ is large, then~$h$ will be minimized quickly and satisfaction of~$g$ will only be attained afterwards, which prolongs convergence. Conversely, for smaller~$\beta$,
minimizing~$h$ and satisfying~$g$ are weighted more equally, which promotes convergence of both equally.

Varying the communication rate has a significant impact on the number of iterations required 
and the behavior of error in agents' updates. Communication rate can even outweigh the benefits of a large~$\beta$, shown in Figure~2. 
This reveals that faster convergence can be achieved by both improving communication or increasing the diagonal dominance of the problem. 
However, if communication is poor, increasing diagonal dominance may not improve results by much, which suggests
that even favorable problem structure does not eliminate the impact of asynchrony. 
Oscillations in the proximity of the optimum are also amplified by large values of~$\beta$.

\section{Conclusion} \label{sec:concl}
After exploring a counterexample to asynchronous dual communications, Algorithm~\ref{alg:2} presents an asynchronous primal-dual approach that is asynchronous in primal updates and communications and asynchronous in distributed dual updates. A numerical example illustrates the effect diagonal dominance has with other parameters upon convergence. Future work will apply the algorithm to large-scale machine learning problems and explore reducing the asynchrony penalty.

\appendix
\emph{Proof of Theorem 1:} Consider the quadratic program
\begin{align}
\textnormal{minimize } &\frac{1}{2}x^TQx + r^Tx \\
\textnormal{subject to } &Ax \leq b, \quad x \in X,
\end{align}
where~$\R^{n \times n} \ni Q = Q^T \succ 0$, $r \in \R^n$, $A \in \R^{m \times n}$,
and~$b \in \R^m$. We take~$X$ sufficiently large in a sense
to be described below. 

Because~$Q$ is symmetric and positive definite, its eigenvectors can be orthonormalized,
and we denote these eigenvectors by~$v_1, \ldots, v_n$. To construct an example,
suppose that~$A \in \R^{m \times n}$ has row~$i$ equal to the~$i^{th}$
normalized eigenvector of~$Q$. 
To compute~$\hat{x}_1 := \argmin_{x \in X} L_{\delta}(x,, \mu^1)$
we set~$\nabla_{x} L_{\delta}(\hat{x}_1, \mu^1) = 0$. 
Expanding and solving, we find~${\hat{x}_1 = -Q^{-1}A^T\mu^1}$,
where we assume that~$X$ is large enough to contain this point. 
Similarly, we find~${\hat{x}_2 = -Q^{-1}A^T\mu^2}$,
where we also assume~$\hat{x}_2 \in X$. 
Then 
\begin{equation}
\|\hat{x}_1 \!-\! \hat{x}_2\|_2 \!=\! \|Q^{-1}\!A^T(\mu^2 \!\!-\!\! \mu^1)\|_2 \!\geq\! \sigma_{min}\big(Q^{-1}\!A^T\big) \|\mu^1 \!-\! \mu^2\|_2, \label{eq:examplemainineq}
\end{equation}
where~$\sigma_{min}(\cdot)$ is the minimum singular value. Expanding, 
\begin{equation}
\sigma_i^2\big(Q^{-1}A^T\big) = \lambda_i\big(AQ^{-T}Q^{-1}A^T\big) = \lambda_i\big(AQ^{-2}A^T\big), 
\end{equation}
and
$AQ^{-2}A^T = \textnormal{diag}\left(\frac{1}{\lambda_1(Q)^2}, \ldots, \frac{1}{\lambda_n(Q)^2}\right)$,
which follows from the fact that we have orthonormalized~$Q$'s eigenvectors. 
Then~$\sigma_{min}\big(Q^{-1}A^T\big) = \frac{1}{\lambda_{max}(Q)}$. 
Using this above, we find 
\begin{equation}
\|\hat{x}_1 - \hat{x}_2\|_2 \geq \frac{1}{\lambda_{max}(Q)} \|\mu^1 - \mu^2\|_2.
\end{equation}
To enforce~$\|\hat{x}_1 - \hat{x}_2\| > L$, 
we ensure that~$\frac{\epsilon}{\lambda_{max}(Q)} > L$,
which is attained for any matrix satisfying~$\lambda_{max}(Q) < \frac{\epsilon}{L}$. 
\hfill $\blacksquare$

\bibliographystyle{plain}
\bibliography{sources2}

\end{document}